\documentclass[letterpaper,10 pt,conference]{IEEEtran}
\IEEEoverridecommandlockouts
\usepackage{cite}
\usepackage{amsmath,amssymb,amsfonts}
\usepackage{algorithmic}
\usepackage{graphicx}
\usepackage{textcomp}
\usepackage{xcolor}
\usepackage{siunitx}
\usepackage{longtable,tabularx}
\usepackage{mathrsfs} 
\usepackage{subcaption}
\usepackage{tikz}
\usepackage[colorlinks=true, allcolors=blue]{hyperref}

\def\BibTeX{{\rm B\kern-.05em{\sc i\kern-.025em b}\kern-.08em
    T\kern-.1667em\lower.7ex\hbox{E}\kern-.125emX}}
\begin{document}

\title{Neural-Network-Based Optimal Guidance for Lunar Vertical Landing\\
\thanks{This research was supported by the National Natural Science Foundation of China under grant No. 62088101.} 
\thanks{The authors are with the School of Aeronautics and Astronautics, Zhejiang University, Hangzhou, China.}
\thanks{ Corresponding author: Zheng Chen (E-mail: z-chen@zju.edu.cn)}
}

\author{\IEEEauthorblockN{
Kun Wang, Zheng Chen*, Fangmin Lu, and Jun Li}
}

\maketitle

\begin{abstract}
This paper addresses an optimal guidance problem concerning the vertical landing of a lunar lander with the objective of minimizing fuel consumption. The vertical landing imposes a final attitude constraint, which is treated as a final control constraint. To handle this constraint, we propose a nonnegative small regularization term to augment the original cost functional. This ensures the satisfaction of the final control constraint in accordance with Pontryagin's Minimum Principle. By leveraging the necessary conditions for optimality, we establish a parameterized system that facilitates the generation of numerous optimal trajectories, which contain the nonlinear mapping from the flight state to the optimal guidance command. Subsequently, a neural network is trained to approximate such mapping. Finally, numerical examples are presented to validate the proposed method.
\end{abstract}


\section{Introduction}
\label{intro}
Exploration of the Moon has been receiving a great deal of attention in recent years. To deploy the lunar lander onto the lunar surface, soft landing of the lunar lander is one of the most fundamental technologies. The soft landing typically involves two phases: de-orbit burn phase and powered descent phase \cite{leeghim2016optimal}. The powered descent phase aims to use the retro-propulsion determined by the guidance command to meet the required final condition \cite{lu2018propellant}. Given the importance of the powered descent phase, considerable efforts have been dedicated to designing the Powered Descent Guidance (PDG).

During the Apollo era, when computational capabilities were quite limited, the landing problem was often simplified to derive an analytical form of the PDG for onboard implementation \cite{cherry1964general}. However, one major drawback of the analytical PDG is the suboptimal fuel consumption. Therefore, there has been growing interest in recent years in deriving the Fuel-Optimal PDG (FOPDG). The FOPDG problem is essentially an Optimal Control Problem (OCP), which can be solved via direct or indirect methods \cite{betts1998survey}. Direct methods reformulate the OCP as a nonlinear programming problem. With the advancements in numerical programming, direct methods, especially convex optimization methods, have shown great potential for onboard implementation in landing problems \cite{accikmecse2013lossless,sagliano2024six}. On the other hand, indirect methods transform the OCP into a Two-Point Boundary-Value Problem (TPBVP) according to the necessary conditions for optimality derived from calculus of variations or Pontryagin’s Minimum Principle (PMP). The resulting TPBVP is then typically solved by the indirect shooting method. Albeit not suitable for onboard implementation, indirect methods are usually used to reveal the structure of the optimal control \cite{lu2018propellant,ito2020throttled,lu2023propellant,wang2024new}. Furthermore, learning-based methods have been employed to generate the real-time optimal solution, where Neural Networks (NNs) are trained to approximate the optimal control; see, e.g., \cite{sanchez2018real,wang2022nonlinear,WANG2024446}. 

However, it is worth noting that most of the methods discussed in the previous paragraph have primarily focused on OCPs without taking into account the final control constraint. From a practical perspective, imposing a final control constraint can ensure smooth system operation \cite{roozegar2018optimal,lee2013polynomial}. In the case of lunar landing, achieving a vertical attitude upon touchdown helps prevent the lunar lander from rolling over and enables the detection of obstructions, facilitating final avoidance maneuvers \cite{mcinnes1995path}. Since the thrust engine is usually fixed with the lunar lander's body, the vertical landing introduces a final attitude constraint, which can be represented by the final steering angle \cite{lu2023propellant}. Consequently, there were some attempts to generate the optimal trajectory ending with a vertical landing \cite{mcinnes1995path,zhou2010optimal,sachan2015fuel,leeghim2016optimal,lu2018propellant}. For instance, in  \cite{leeghim2016optimal}, 
the lunar landing trajectory was divided into two parts, with the final steering angle constraint augmented into the cost functional. As a result, the original TPBVP became a Multi-PBVP. In  \cite{mcinnes1995path}, the concept of gravity-turn was employed to direct the thrust vector opposite to the velocity vector. By integrating the equations of motion for altitude, velocity magnitude, and time analytically as a function of flight path angle, the thrust vector was used as a parameter to shape the trajectory. In  \cite{sachan2015fuel}, to handle the final steering angle constraint, a time-varying matrix was added to the cost functional, and the resulting optimal control problem was solved via  model predictive static programming. 

In this paper, we present a novel approach to handle the final steering angle constraint in lunar vertical landing scenarios. To ensure the satisfaction of this constraint, we propose a modification to the cost functional by introducing a nonnegative small regularization term. Building upon the necessary conditions for optimality derived from PMP, we establish a parameterized system, which allows for generating numerous optimal trajectories. These trajectories capture the nonlinear mapping from the flight state to the optimal steering angle command. Then, an NN is trained to approximate such mapping. Finally, numerical simulations are presented, showing that the proposed method is able to guide the lunar lander to achieve a vertical landing.
\section{Problem Formulation}
\subsection{Dynamical Model and Boundary Conditions}
\begin{figure}[!htp]
\begin{center}
\includegraphics[scale=0.25]{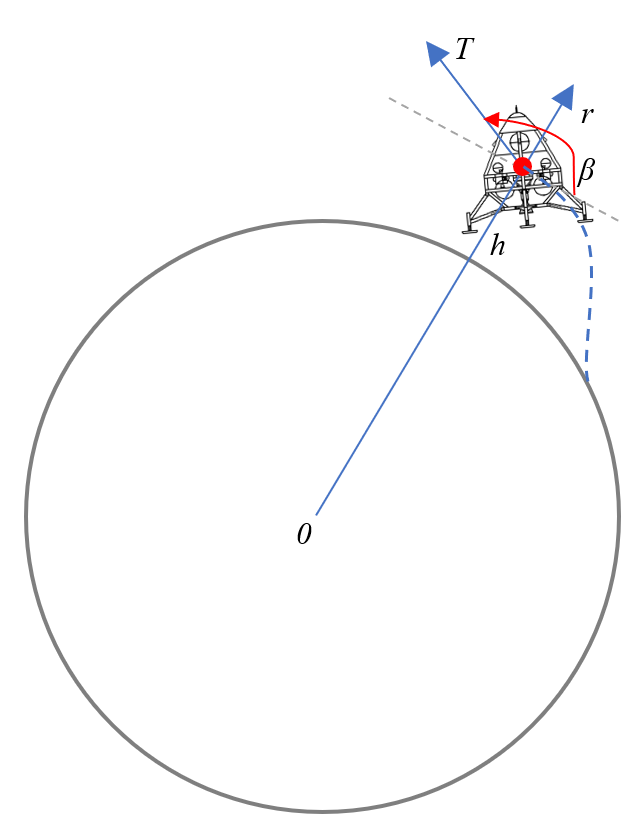}
\caption{Coordinate system for lunar vertical landing.}\label{Fig:frame}
\end{center}
\end{figure}

Consider the planar motion of a lunar lander, as depicted in Fig.~\ref{Fig:frame}. 
Here, it is assumed that the Moon is a regular spherical body and we do not consider the Moon's rotation. The point $O$ is fixed at the center of the Moon. Denote by $r \ge R_0$ ($R_0$ is the radius of the Moon) the radial distance between $O$ and the lunar lander. Consequently, the altitude of the lunar lander above the lunar surface can be expressed as $h = r - R_0$. The lunar lander is propelled by an engine with a constant thrust magnitude $T$ \cite{leeghim2016optimal,pengkun2019}. The thrust steering angle, denoted by $\beta \in [0, 2\pi]$, represents the angle between the local horizontal line of the lunar lander and the thrust vector. It is important to note that since the landing site is not fixed, the optimal trajectory can be rotated. Then, the dynamics of the lunar lander is governed by \cite{leeghim2016optimal}
\begin{align}
\begin{cases}
\dot{r}(t) =  v(t),\\
\dot{u}(t) = -\frac{u(t)v(t)}{r(t)} + \frac{T}{m(t)} \cos \beta(t), \\
\dot{v}(t) =  \frac{u^2(t)}{r(t)}-\frac{\mu}{r^2(t)} + \frac{T}{m(t)} \sin \beta(t),\\
\dot{m}(t) = -\frac{T}{I_{sp}g_e},
\end{cases}
\label{LunarLander:DiffEqution}
\end{align}
where $t \ge 0$ represents time, $u$ the transverse speed and $v$ the radial speed. $\mu$ is the gravitational constant of the Moon, and $m$ is the mass of the lunar lander. $I_{sp}$ denotes the specific impulse of the lunar lander's propeller, while $g_e$ represents the Earth's gravitational acceleration at sea level. The initial condition of the lunar lander is
\begin{align}
r(0) = r_0, ~u(0) = u_0, ~v(0) = v_0, ~m(0) = m_0.
\label{InitialState}
\end{align}

From Fig.~\ref{Fig:frame}, we have $r(t_f) = R_0$ at touchdown ($t_f$ is the free final time). We expect to have zero speed at touchdown, i.e., $u(t_f) = 0$ and $v(t_f) = 0$. Therefore, the final condition can be given by
\begin{align}
r(t_f) = R_0, ~u(t_f) = 0, ~v(t_f) = 0.
\label{FinalState}
\end{align}

For better numerical conditioning, constants $R_0$, $\sqrt\frac{\mu}{R_0}$, $\sqrt\frac{\mu}{R_0}$, $m_0$, $\sqrt\frac{{R_0}^3}{\mu}$, and $\frac{m_0 \mu}{{R_0}^2}$ are employed to normalize $r$, $u$, $v$, $m$, $t$, and $T$, respectively. To avoid notation abuse, we still use the same notation as in Eq.~(\ref{LunarLander:DiffEqution}) for the dimensionless counterpart in the rest of the paper.
\subsection{Optimal Guidance Problem Formulation} 
Because the steering angle defines the body attitude of the lunar lander \cite{lu2023propellant}, to ensure a vertical landing, the final steering angle should satisfy the following constraint:
\begin{align}
\beta(t_f) = \frac{\pi}{2}.
\label{EQ:steeringcondition}
\end{align}
The cost functional to be minimized is
\begin{align}
J = \int_0^{t_f} 1 ~\mathrm{d}t,
\label{EQ:cost}
\end{align}
which is equivalent to minimizing the fuel consumption since the thrust magnitude is constant.

The optimal guidance problem amounts to finding the real-time steering angle $\beta(t)$ and the final time $t_f$, to steer the dynamical system specified by Eq.~(\ref{LunarLander:DiffEqution}) from the initial condition in Eq.~(\ref{InitialState}) to the final condition in Eq.~(\ref{FinalState}), such that the cost functional in Eq.~(\ref{EQ:cost}) is minimized while the constraint in Eq.~(\ref{EQ:steeringcondition}) is met.
It is important to note that Eq.~(\ref{EQ:steeringcondition}) poses a final control constraint, which can be challenging to handle using conventional optimal control theory. In the following, we will propose an easy-to-implement method to address this constraint. By leveraging the necessary conditions for optimality derived from PMP, we establish a parameterized system, as previously done in \cite{wang2022nonlinear,WANG2024446}. This enables one to readily obtain numerous optimal trajectories. Then, an NN is trained to approximate the optimal steering angle command.
\section{Parameterization of Optimal Trajectories}\label{SE:properties}
\subsection{Handling the Final Steering Angle Constraint}
We introduce a regularization term, defined as follows:
\begin{align}
\Delta (r,\beta,\delta,\epsilon) :=  \frac{1}{2}\delta\exp(1-r)\frac{(\beta-\frac{\pi}{2})^2}{r-1+\epsilon},
\label{EQ:new_cost_term}
\end{align}
where $\delta$ is a small positive constant and $\epsilon$ is a very small positive constant to avoid numerical singularity. By augmenting this regularization term to the cost functional in Eq.~(\ref{EQ:cost}), the augmented cost functional $\tilde{J}$ becomes
\begin{align}
\begin{split}
&\tilde{J} = \int_0^{{t}_f} (1+\Delta) ~\mathrm{d}{t} = \\
&\int_0^{{t}_f} \left\{ 1 + \frac{1}{2}\delta\exp[1-r(t)]\frac{[\beta(t)-\frac{\pi}{2}]^2}{r(t)-1+\epsilon}\right\}  ~\mathrm{d}{t}. 
\end{split}
\label{EQ:new_cost}
\end{align}

Next, we will demonstrate a property that holds for the regularization term $\Delta (r,\beta,\delta,\epsilon)$.

{\it Property:} The term $\Delta (r,\beta,\delta,\epsilon)$ remains nonnegative and sufficiently small for ${t} \in [0, {t}_f]$.

{\it Proof:} Before touchdown, we have ${r}({t}) > 1$, and since $\delta$ is a small positive constant, $\delta \exp[1-r(t)]$ becomes a very small positive number. Additionally, it is evident that $\frac{[\beta(t)-\frac{\pi}{2}]^2}{r(t)-1+\epsilon} \geq 0$ holds for the entire landing process.
To minimize the cost functional $\tilde {J}$, the nonnegative term $\frac{[\beta(t)-\frac{\pi}{2}]^2}{r(t)-1+\epsilon}$ cannot approach infinity.
At touchdown ${t} = {t}_f$, we have ${r}({t})=1$. In this case,  both $\lim\limits_{{t}\to {t}_f} {\beta}({t}) = \frac{\pi}{2}$ and  $\lim\limits_{{t}\to {t}_f} \frac{[\beta(t)-\frac{\pi}{2}]^2}{r(t)-1+\epsilon} = 0$ have to hold to ensure that the regularization term $\Delta (r,\beta,\delta,\epsilon)$ does not become extremely large. This completes the proof. $\square$

With this property, the final steering angle constraint in Eq.~(\ref{EQ:steeringcondition}) will be automatically satisfied without imposing a significant penalty on the original cost functional in Eq.~(\ref{EQ:cost}). In the following sections, we will derive the necessary conditions for optimality w.r.t. the augmented cost functional.
\subsection{Necessary Conditions for Optimality}
Rewrite the dynamics in Eq.~(\ref{LunarLander:DiffEqution}) as
\begin{align}
\dot{{\boldsymbol x}}( t) = \boldsymbol f({\boldsymbol x}, {\beta}, t),
\label{EQ:new_dynamics}
\end{align}
where $\boldsymbol x = [r,u,v,m]^T$ and $\boldsymbol f:\mathbb{R}^4 \times \mathbb{R} \times \mathbb{R}_0^+ \rightarrow \mathbb{R}^4$ is the dynamics defined in Eq.~(\ref{LunarLander:DiffEqution}).

Denote by $\boldsymbol {p}_x = [p_r,p_u,p_v,p_m]^T$ the co-state vector.
Then, the Hamiltonian is built as
\begin{align}
\begin{split}
&\mathscr{H} = \boldsymbol {p}^T_x \boldsymbol f + 1 + \Delta = p_r v + p_u
(-\frac{uv}{r} + \frac{T}{m} \cos \beta) + \\
 &p_v(
\frac{u^2}{r}-\frac{\mu}{r^2} +  \frac{T}{m} \sin \beta) - p_m \frac{T}{I_{sp}g_e} + \\
& [ 1 + \frac{1}{2}\delta\exp(1-r)\frac{(\beta-\frac{\pi}{2})^2}{r-1+\epsilon}]
. 
\end{split}
\label{EQ:Ham_modified}
\end{align}
Based on PMP \cite{Pontryagin}, the dynamics of the co-state vector is  
\begin{align*}
\dot{\boldsymbol {p}}_x(t) = -\frac{\mathscr{H}}{\partial \boldsymbol {x}(t)},
\end{align*}
that is 
\begin{align}
\begin{cases}
\dot{{p}}_r( t) =  -\frac{2\mu p_v}{r^3} + \frac{p_v u^2 - p_u u v}{r^2} 
+ \frac{1}{2}\delta\exp(1-r)\frac{(\beta-\frac{\pi}{2})^2(r+\epsilon)}{(r-1+\epsilon)^2},\\
\dot{{p}}_u( t) = \frac{p_u v - 2 p_v u}{r},\\
\dot{{p}}_v( t) =  \frac{up_u}{r} - p_r,\\
\dot{p}_m( t) =  \frac{T}{{m}^2}(p_u \cos \beta + p_v \sin \beta).
\end{cases}
\label{EQ:dot_p_modified}
\end{align}
The optimal steering angle satisfies that
\begin{align}
\frac{\partial {\mathscr{H}}}{\partial {{\beta}}( t)} = 0.
\label{EQ:dH/dtheta_modi}
\end{align}
Explicitly rewriting Eq.~(\ref{EQ:dH/dtheta_modi}) yields
\begin{align}
\begin{split}
 &\frac{T}{{m}(t)}[-p_u(t)\sin\beta(t) + p_v(t) \cos \beta(t)] + \\
 &\delta \exp[1-r(t)]\frac{\beta(t)-\frac{\pi}{2}}{r(t)-1+\epsilon} = 0.
 \end{split}
\label{EQ:dH/dtheta_modi_real}
\end{align}

It is worth noting that if we disregard the final steering angle constraint (i.e., $\delta=0$), the optimal steering angle can be analytically determined in terms of $p_u$ and $p_v$. However, when $\delta > 0$, the optimal steering angle needs to be numerically solved using Eq.~(\ref{EQ:dH/dtheta_modi_real}). This equation is a transcendental equation that may have multiple zeros. Thus, conventional iterative methods like Newton's method or the bisection method may not find the desired zero. Based on the approach presented in \cite{zheng2021time}, we shall show how to find the right zero for Eq.~(\ref{EQ:dH/dtheta_modi_real}).

Differentiating Eq.~(\ref{EQ:dH/dtheta_modi_real}) w.r.t.  $ \beta$ leads to
\begin{align}
  \frac{T}{{m}(t)}[-p_u(t)\cos\beta(t) - p_v(t) \sin \beta(t)] + \frac{\delta \exp[1-r(t)]}{r(t)-1+\epsilon} = 0.
\label{EQ:dH/dtheta_modi_real_zeros}
\end{align}
By substituting the half-angle formulas
\begin{align*}
\sin{{\beta}} = \frac{2\tan\frac{{{\beta}}}{2}}{1+\tan^2\frac{{{\beta}}}{2}} ~~ {\rm and} ~~
\cos{{\beta}} = \frac{1-\tan^2\frac{{{\beta}}}{2}}{1+\tan^2\frac{{{\beta}}}{2}},
\end{align*}
into Eq.~(\ref{EQ:dH/dtheta_modi_real_zeros}), we have that $\tan\frac{{{\beta}}}{2}$ is a zero of the following quadratic polynomial in terms of $x$:
\begin{align}
\begin{split}
&\left\{\frac{\delta \exp[1-r(t)]}{r(t)-1+\epsilon} + p_u(t)\frac{T}{{m}(t)}\right\}x^2-2{p}_{v}(t)\frac{T}{{m}(t)}x + \\
&\left\{\frac{\delta \exp[1-r(t)]}{r(t)-1+\epsilon} - p_u(t)\frac{T}{{m}(t)}\right\} = 0.
\end{split}
\label{EQ:dH/dtheta_modi_real_quad}
\end{align}

Note that the two roots of Eq.~(\ref{EQ:dH/dtheta_modi_real_quad}) can be obtained either by radical expressions or standard polynomial solvers. Then a simple bisection method can be applied to find all the zeros of Eq.~(\ref{EQ:dH/dtheta_modi_real}) by comparing the values for Eq.~(\ref{EQ:dH/dtheta_modi_real}) when ${{\beta}}$ takes values of $0$, $2\pi$, and the two roots of Eq.~(\ref{EQ:dH/dtheta_modi_real_quad}). After finding all the zeros of Eq.~(\ref{EQ:dH/dtheta_modi_real}), the optimal value for ${{\beta}}$ is determined by selecting the one that minimizes the Hamiltonian in Eq.~(\ref{EQ:Ham_modified}).

The free final mass leads to the transversality condition:
\begin{align}
{p}_{m}({t}_f) = 0.
\label{EQ:Transversality1}
\end{align}
In addition, $\mathscr H$ does not explicitly depend on time and the final time is free, leading to the following stationary condition:
\begin{align}
\mathscr H(t_f) = 0.
\label{EQ:stat1}
\end{align}

Eqs.~(\ref{LunarLander:DiffEqution}), (\ref{EQ:dot_p_modified}), (\ref{EQ:dH/dtheta_modi_real}), (\ref{EQ:Transversality1}) and (\ref{EQ:stat1}) constitute a set of ordinary differential equations in terms of the states and co-states. The resulting TPBVP is formulated as:
\begin{align}
\boldsymbol \psi(\boldsymbol{p}_{x_0}, t_f) = [r(t_f)-R_0,u(t_f),v(t_f),p_m(t_f),\mathscr H(t_f)],
\label{EQ:TPBVP}
\end{align}
where $\boldsymbol \psi$ is the shooting function, while $\boldsymbol{p}_{x_0}$ and $t_f$ is the initial guess of co-state vector $\boldsymbol{p}_{x}$ and the final time, respectively.

However, solving Eq.~(\ref{EQ:TPBVP}) can be time-consuming and may encounter convergence issues due to the lack of physical significance in $\boldsymbol{p}_{x_0}$ \cite{chen2019nonlinear}. In the next subsection, we will introduce a parameterized system that allows for obtaining optimal trajectories efficiently, enabling the generation of numerous optimal trajectories required for training the NN.
\subsection{Parameterized System} 
Define a new variable $\tau$ as
\begin{align}
\tau = t_f - t, t \in [0,t_f].
\label{Eq:tao}
\end{align}
We establish a first-order ordinary differential system 
\begin{align}
\begin{cases}
\dot{\bar {\boldsymbol x}} = -\boldsymbol f(\boldsymbol {\bar x}, {\bar \beta}, \tau), \\
\dot{\bar{\boldsymbol{p}}}_{x} = \frac{\partial {\bar {\mathscr{H}}}(\tau)}{\partial \boldsymbol {\bar x}(\tau)},
\label{EQ:pareEquation}
\end{cases}
\end{align}
where $\boldsymbol {\bar x} = [\bar r,\bar u,\bar v, \bar m]^\top$, $\bar{\boldsymbol{p}}_{x} =  [\bar{p}_r,\bar{p}_u,\bar {p}_v, \bar {p}_m]^\top$,  and $\bar {\mathscr{H}} = \bar{\boldsymbol {p}}^T_x \boldsymbol f + 1 + \bar{\Delta}$. 
Meanwhile, $\bar \beta$ satisfies
\begin{align}
\begin{split}
&\frac{T}{\bar {m}(\tau)}[-\bar{p}_u(\tau)\sin\bar{\beta}(\tau) + \bar{p}_v(\tau) \cos \bar{\beta}(\tau)] + \\
&\delta \exp[1-\bar{r}(\tau)]\frac{\bar{\beta}(\tau)-\frac{\pi}{2}}{\bar{r}(\tau)-1+\epsilon} = 0,
\end{split}
\label{EQ:optimal_psi_para}
\end{align}

The initial condition at $\tau = 0$ for the system in Eq.~(\ref{EQ:pareEquation}) is defined as
\begin{align}
\begin{cases}
\bar r(0) = 1, ~\bar u(0) = 0, ~\bar v (0) = 0, ~\bar m(0) = \bar {m}_0, \\
\bar{p}_r(0) = \bar{p}_{r_0},~\bar{p}_u(0) = \bar{p}_{u_0},~\bar{p}_v(0) = \bar{p}_{v_0}, ~\bar{p}_m(0) = 0,
\label{EQ:Initial_pareEquation}
\end{cases}
\end{align}
where the triple $(\bar {p}_{r_0}, \bar {p}_{u_0}, \bar {p}_{v_0})$ is arbitrary. The value for $\bar {m}_0$ is chosen to meet 
\begin{align}
\bar {\mathscr{H}}(0) = 0.
\label{EQ:Hpara}
\end{align}
Substituting Eq.~(\ref{EQ:Initial_pareEquation}) into Eq.~(\ref{EQ:Hpara}) leads to
\begin{align}
\begin{split}
\bar {\mathscr{H}}(0) = &\frac{T}{\bar {m}_0}[\bar {p}_{u_0} \cos \bar{\beta}(0) + \bar {p}_{v_0} \sin \bar{\beta}(0)] -\mu \bar {p}_{v_0} \\
&+ 1 + \frac{1}{2}\delta\frac{[\bar{\beta}(0)-\frac{\pi}{2}]^2}{\epsilon}. 
\end{split}
\label{EQ:Ham_para_initial}
\end{align}

Therefore, for any pair $(\bar {p}_{u_0}, \bar {p}_{v_0})$, the values for $\bar {m}_0$  and $\bar{\beta}(0)$ can be determined by solving Eqs.~(\ref{EQ:optimal_psi_para}) and (\ref{EQ:Ham_para_initial}). Consequently, with an arbitrary $\bar {p}_{r_0}$, an optimal trajectory can be readily obtained by propagating the system in Eq.~(\ref{EQ:pareEquation}) with the initial condition specified in Eq.~(\ref{EQ:Initial_pareEquation}). Let the pair
$(\boldsymbol {\bar x}(\tau,\bar {p}_{r_0}, \bar {p}_{u_0}, \bar {p}_{v_0}),\bar{\boldsymbol{p}}_{x}(\tau,\bar {p}_{r_0}, \bar {p}_{u_0}, \bar {p}_{v_0})) \in \mathbb{R}^{8}$
for $\tau \in [0,t_f]$ denote the solution to the parameterized system in Eq.~(\ref{EQ:pareEquation}) with the initial condition specified in Eq.~(\ref{EQ:Initial_pareEquation}). In this case, it is evident that this pair meets all the necessary conditions for optimality. Therefore, an optimal trajectory can be obtained by arbitrarily choosing a triple $(\bar {p}_{r_0}, \bar {p}_{u_0}, \bar {p}_{v_0})$ and propagating the parameterized system in Eq.~(\ref{EQ:pareEquation}) with the initial condition in Eq.~(\ref{EQ:Initial_pareEquation}). 
Moreover, the optimal steering angle command $\bar{\beta}$ can be determined by Eq.~(\ref{EQ:dH/dtheta_modi_real}). Denote by $f_{\bar\beta}$ the nonlinear mapping $\boldsymbol{\bar x} \longmapsto \bar \beta$. Based on the universal approximation theorem \cite{HORNIK1989359}, if we have a large number of sampled data representing the relationship $f_{\bar\beta}$, a well-trained NN will be capable of representing this relationship in an accurate way. This will be further elucidated in the subsequent section.
\section{Neural-Network-Based Real-Time Solution}
To train an NN for generating the optimal steering angle command, it is necessary to create a training dataset. For this purpose, a nominal trajectory is required. 
The radius of the Moon is $R_0 = 1,738$ km. The propulsion system of the lunar lander is characterized by $I_{sp} = 300$ s and $T = 1,500$ N. $g_e$ is equal to 9.81 $\rm m/s^2$.  The constants in the regularization term defined in Eq.~(\ref{EQ:new_cost_term}) are $\delta = 1.0 \times 10^{-5}$ and $\epsilon = 1.0 \times 10^{-8}$. The nominal initial condition is set as: $r_0 = 1,753$ km, $u_0 = 1,6795 $ m/s, $v_0 = 0$ m/s, $m_0 = 600$ kg.

\begin{figure}[!htp]
\begin{center}
\includegraphics[scale=0.16]{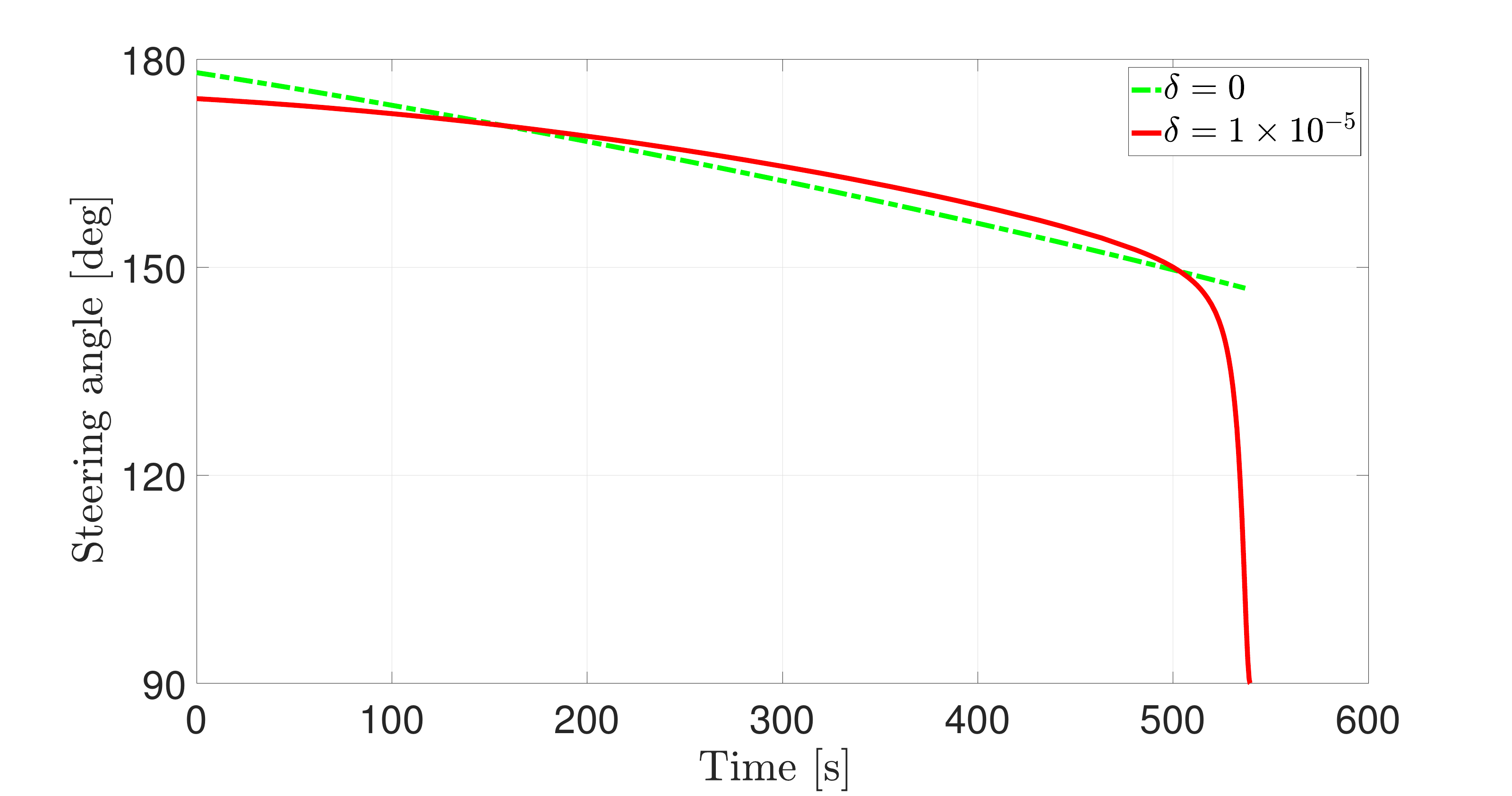}
\caption{Profiles of the steering angle.}\label{Fig:steering_nominal}
\end{center}
\end{figure}

\begin{figure}[!htp]
\begin{center}
\includegraphics[scale=0.16]{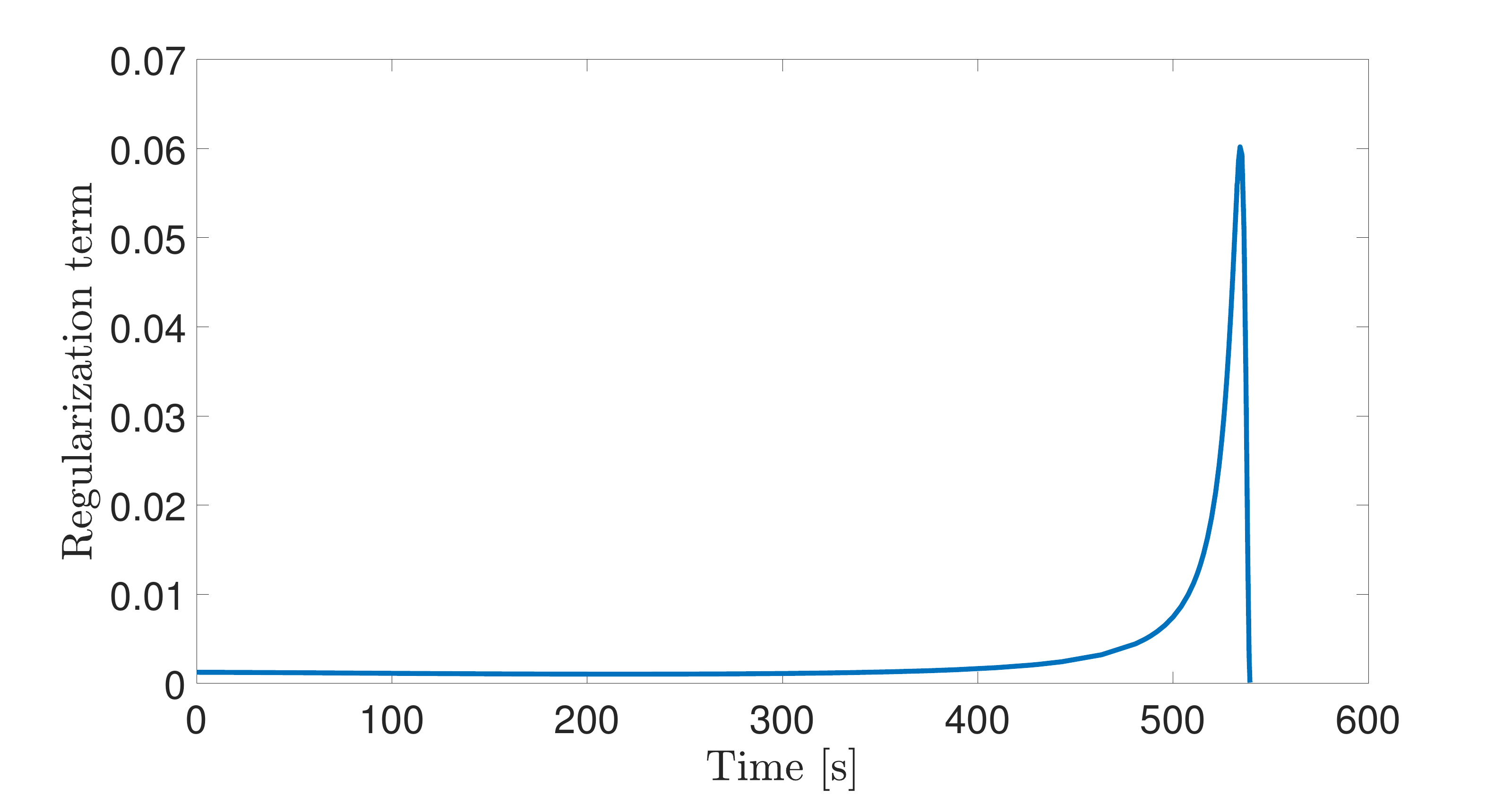}
\caption{Profile of the regularization term.}\label{Fig:term_nominal}
\end{center}
\end{figure}

Initially, we consider the case without the final steering angle constraint, i.e., $\delta=0$. The indirect shooting method is used to solve the corresponding TPBVP defined in Eq.~(\ref{EQ:TPBVP}). A homotopy procedure is adopted to increase $\delta$ from $0$ to $1\times 10^{-5}$. Fig.~\ref{Fig:steering_nominal} illustrates the profiles of the steering angle for $\delta=0$ and $1\times 10^{-5}$. When the final steering angle constraint is not considered, the steering angle exhibits an almost linear behavior. However, for $\delta=1\times 10^{-5}$,  it can be observed that the final steering angle reaches exactly $90$ deg, indicating the fulfillment of the requirement for a vertical landing. Moreover, the final time for $\delta=0$ is $536.90$ seconds, while it increases to $539.29$ seconds for $\delta=1\times 10^{-5}$. The profile of the regularization term $\Delta$ related to $\delta=1\times 10^{-5}$ is depicted in Fig.~\ref{Fig:term_nominal}, indicating that the regularization term remains nonnegative.
As a result, the regularization term incurs an additional fuel consumption of $1.2181$ kg, which is negligible in the context of an initial mass of $600$ kg. The solution related to $\delta=1\times 10^{-5}$ is chosen as the nominal trajectory.

\begin{figure}[!htp]
\begin{center}
\includegraphics[scale=0.16]{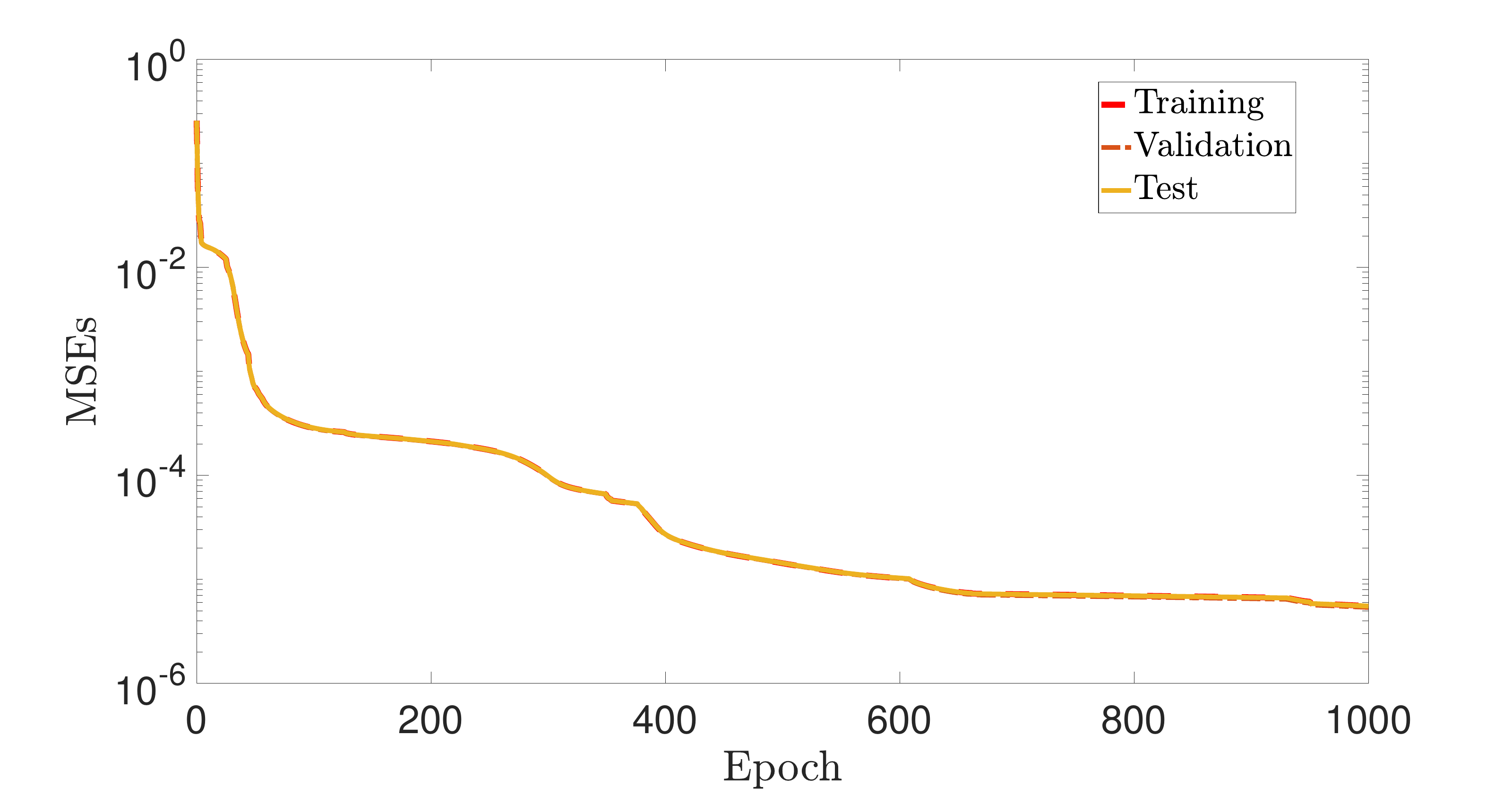}
\caption{Training history of the NN.}\label{Fig:training}
\end{center}
\end{figure}

By applying perturbations over the nominal trajectory, we manage to collect $16,903$ optimal trajectories; from each optimal trajectory, state–steering angle pairs are selected along the trajectory and inserted in the training dataset. An NN with three hidden layers (each of which contains $15$ neurons) is trained to approximate $f_{\bar\beta}$. The sigmoid function is employed for all the hidden layers, and a linear function is used for the output layer. Before training, a split of $70\%$ for training, $15\%$ for validation, and $15\%$ for testing is shuffled for the  dataset samples. The essence of the training is minimizing the loss function, which is quantified as the Mean Squared Error (MSE) between the predicted values from the trained NN and the actual values within the dataset samples. The training is terminated after $1,000$ epochs. Fig.~\ref{Fig:training} illustrates the training progression of the NN.
As a result, the MSEs for the training, validation, and testing sets reduce to $5.39 \times 10^{-6}$. Consequently, the trained NN can not only generate the steering angle command within a constant time, but also offers a closed-form solution, as demonstrated by the numerical simulations in the next section.
\section{Numerical Simulations}
We consider two initial conditions outside the training dataset, as presented in Table~\ref{tab1}. For comparison, the subroutine {\it fsolve} is used to solve the shooting function in Eq.~(\ref{EQ:TPBVP}) for the case of $\delta=0$, which corresponds to the scenario without the final steering angle constraint. On the other hand, for the scenario with the final steering angle constraint (i.e., $\delta=1\times 10^{-5}$), we employ the trained NN to generate the steering angle. In order to ensure a stable output, the simulation using the trained NN is terminated once the altitude drops below $5$ m. It is noteworthy that all the algorithms are implemented on a laptop equipped with an AMD Ryzen 7-5800H CPU operating at 3.2 GHz.

The indirect shooting method takes approximately $1.5100$ seconds to find the open-loop optimal solution for $\delta=0$. It is worth noting that the implicit shooting method also required nearly 1 second to find the solution \cite{pengkun2019}. However, the trained NN is capable of generating the steering angle within $0.0066$ seconds given a flight state. This indicates that the trained NN can provide real-time closed-form solutions.
\begin{table}[htbp]
\caption{Initial Conditions}
\begin{center}
\begin{tabular}{|c|c|c|c|c|}
\hline
Item
 & $r_0$ (km)  & $u_0$ (m/s) & $v_0$ (m/s) & $m_0$ (kg)\\
\hline
Case 1 & $1,9608.12$  & $76.7558$  & $-825.1668$  & $600$  \\
\hline
Case 2 & $1,7407.78$  & $86.3188$  & $-108.7746$  & $376.5833$  \\
\hline
\end{tabular}
\label{tab1}
\end{center}
\end{table}

\begin{figure}[!htp]
\centering
\begin{subfigure}[t]{4.05cm}
\centering
\includegraphics[width = 4.05cm]{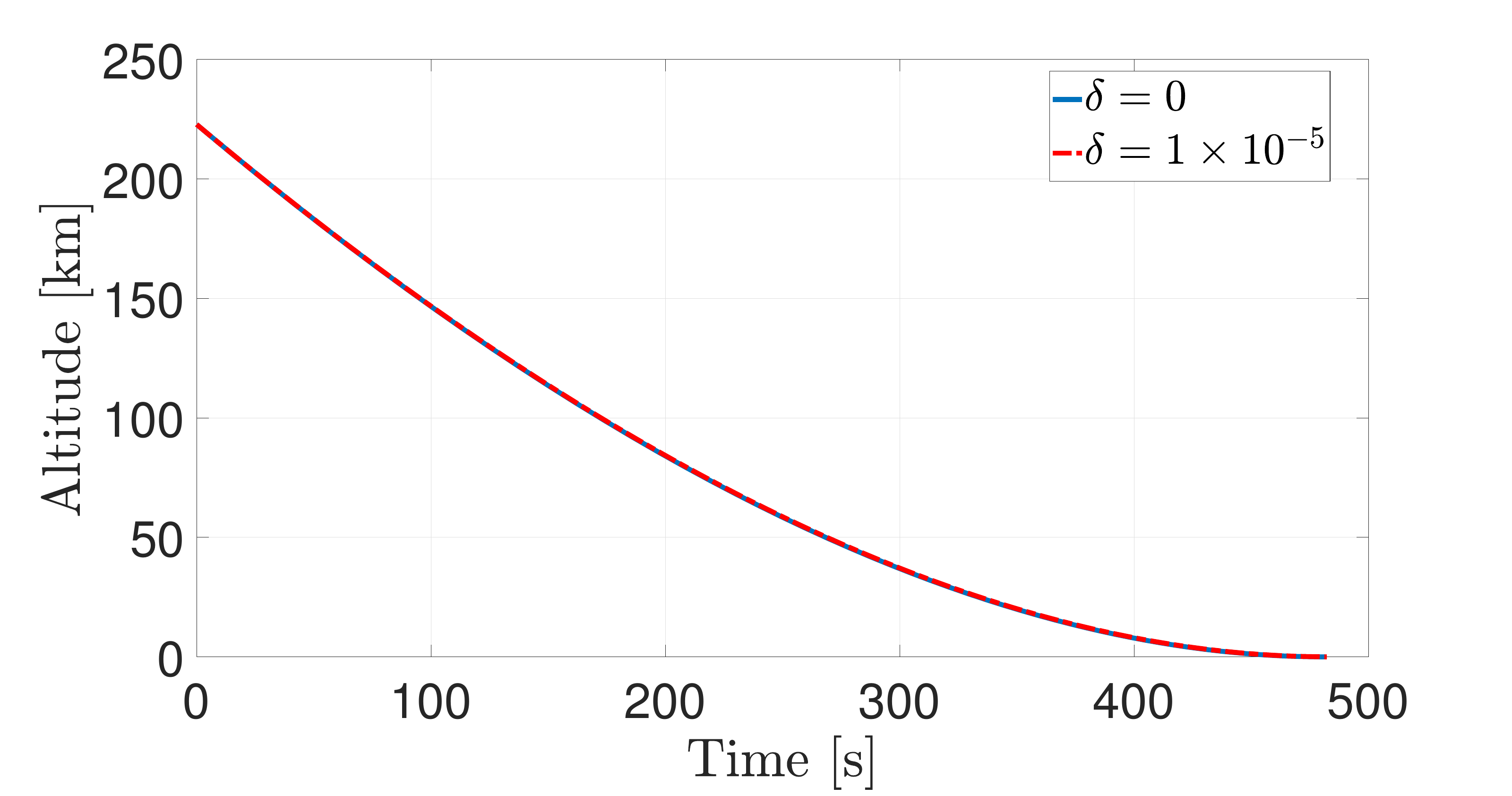}
\caption{Altitude}
\label{Fig:cooperative_control_1_case1}
\end{subfigure}
~~~~~
\begin{subfigure}[t]{4.05cm}
\centering
\includegraphics[width = 4.05cm]{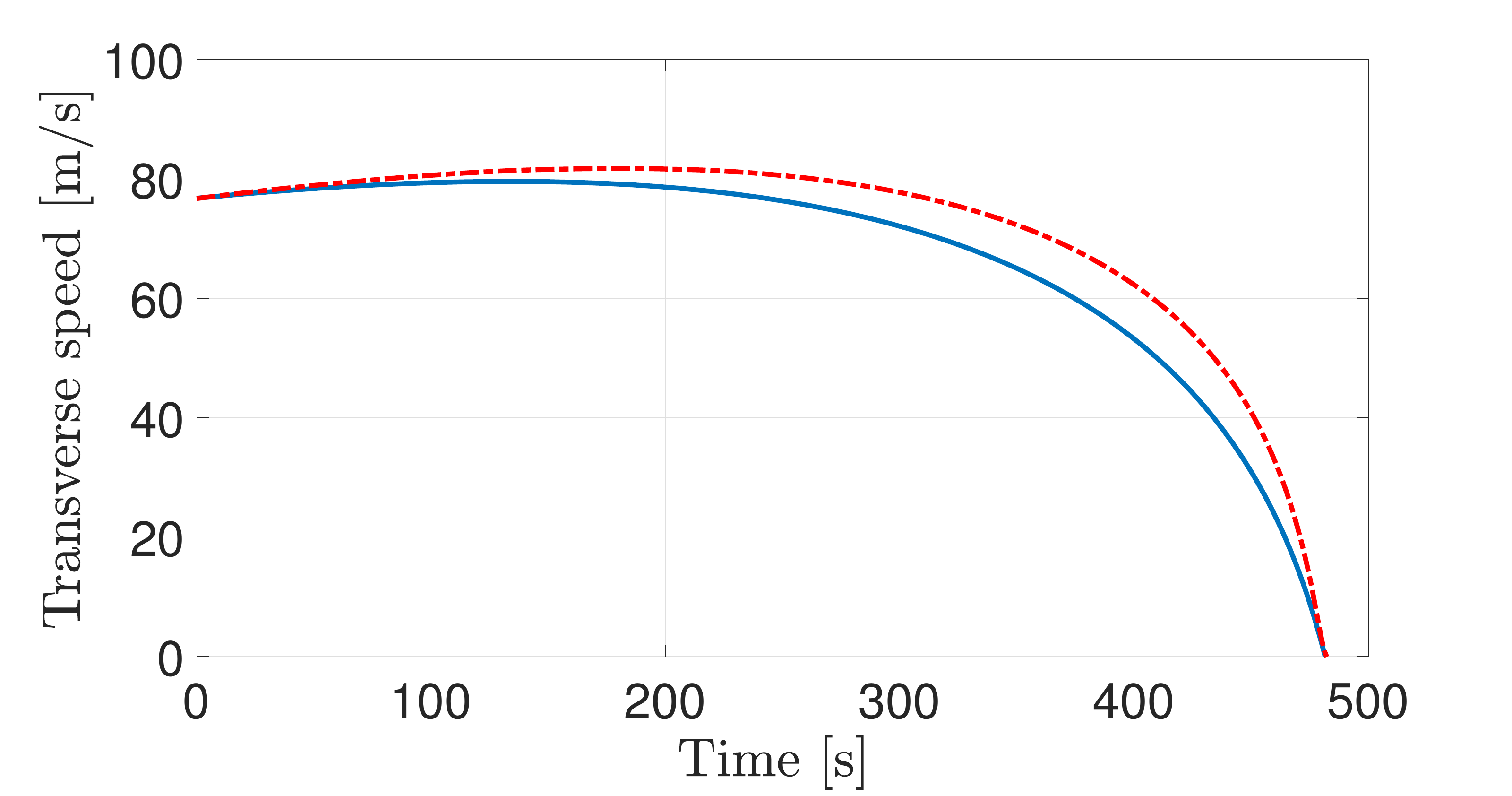}
\caption{Transverse speed}
\label{Fig:cooperative_control_2_case1}
\end{subfigure}\\
\begin{subfigure}[t]{4.05cm}
\centering
\includegraphics[width = 4.05cm]{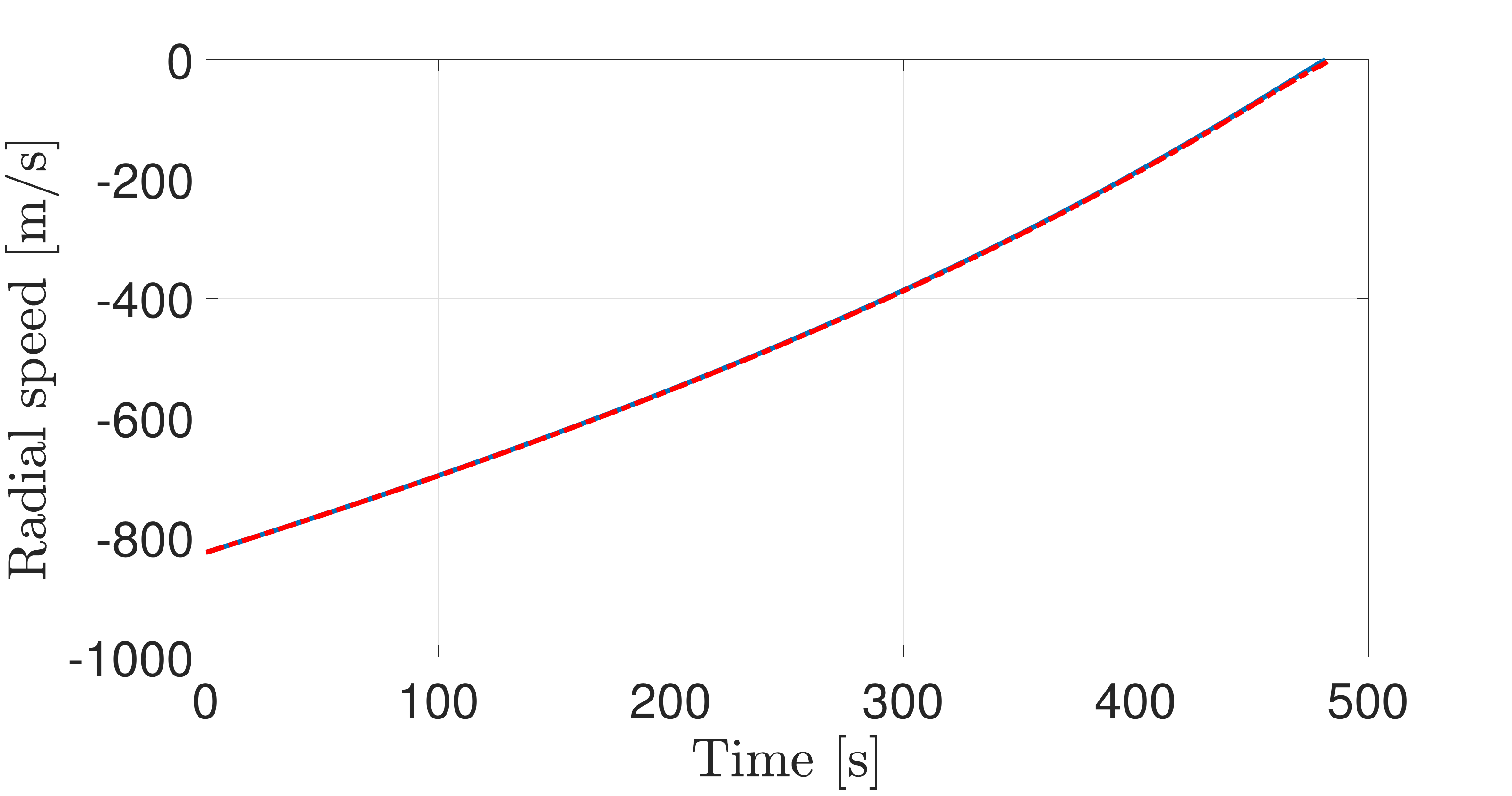}
\caption{Radial speed}
\label{Fig:cooperative_control_3_case1}
\end{subfigure}
~~~~~
\begin{subfigure}[t]{4.05cm}
\centering
\includegraphics[width = 4.05cm]{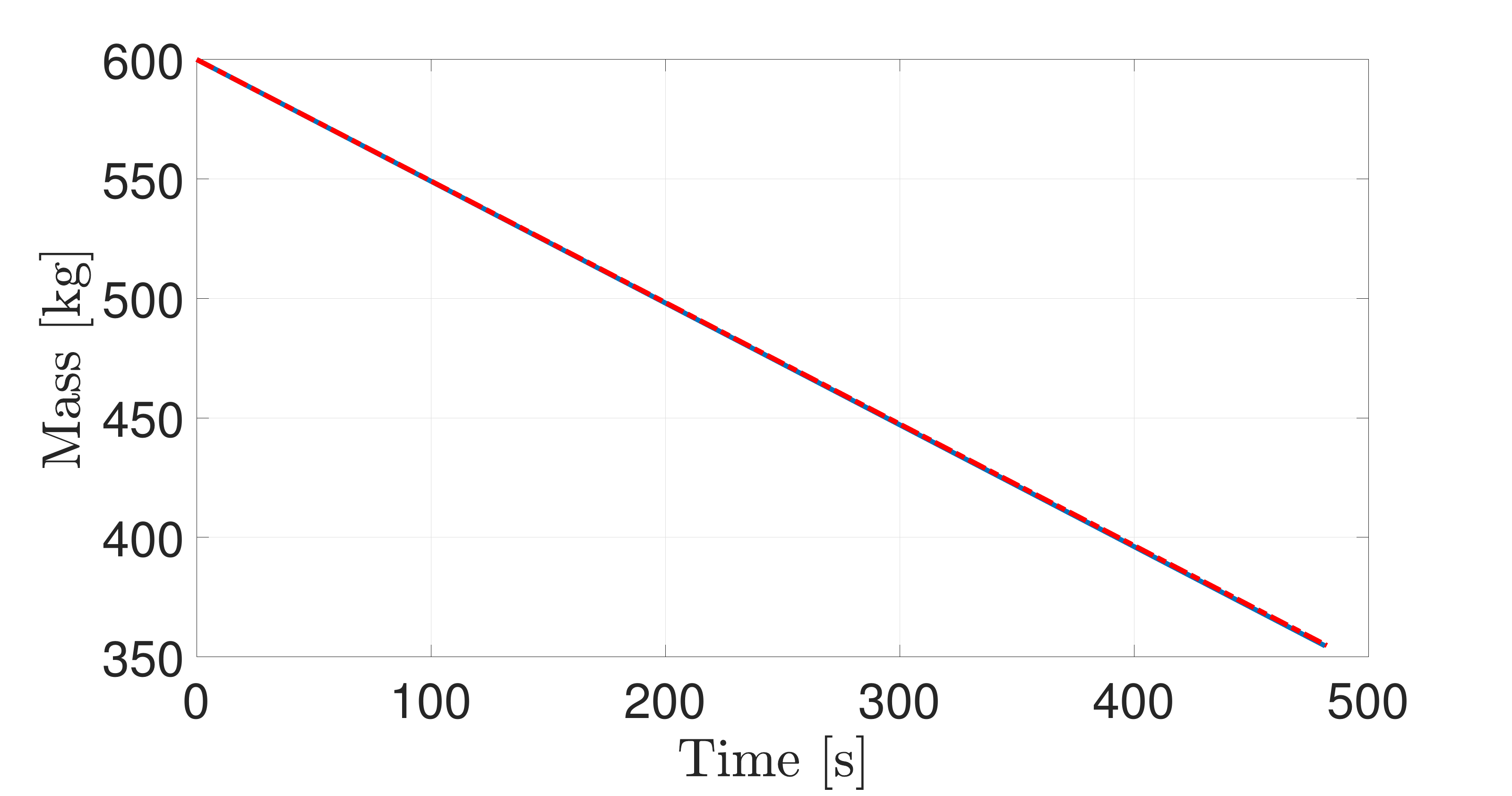}
\caption{Mass}
\label{Fig:cooperative_control_4_case1}
\end{subfigure}
\caption{Profiles of the states without and with the final steering angle constraint for Case 1.}
\label{Fig:cooperative_profile}
\end{figure}

\begin{figure}[!htp]
\begin{center}
\includegraphics[scale=0.16]{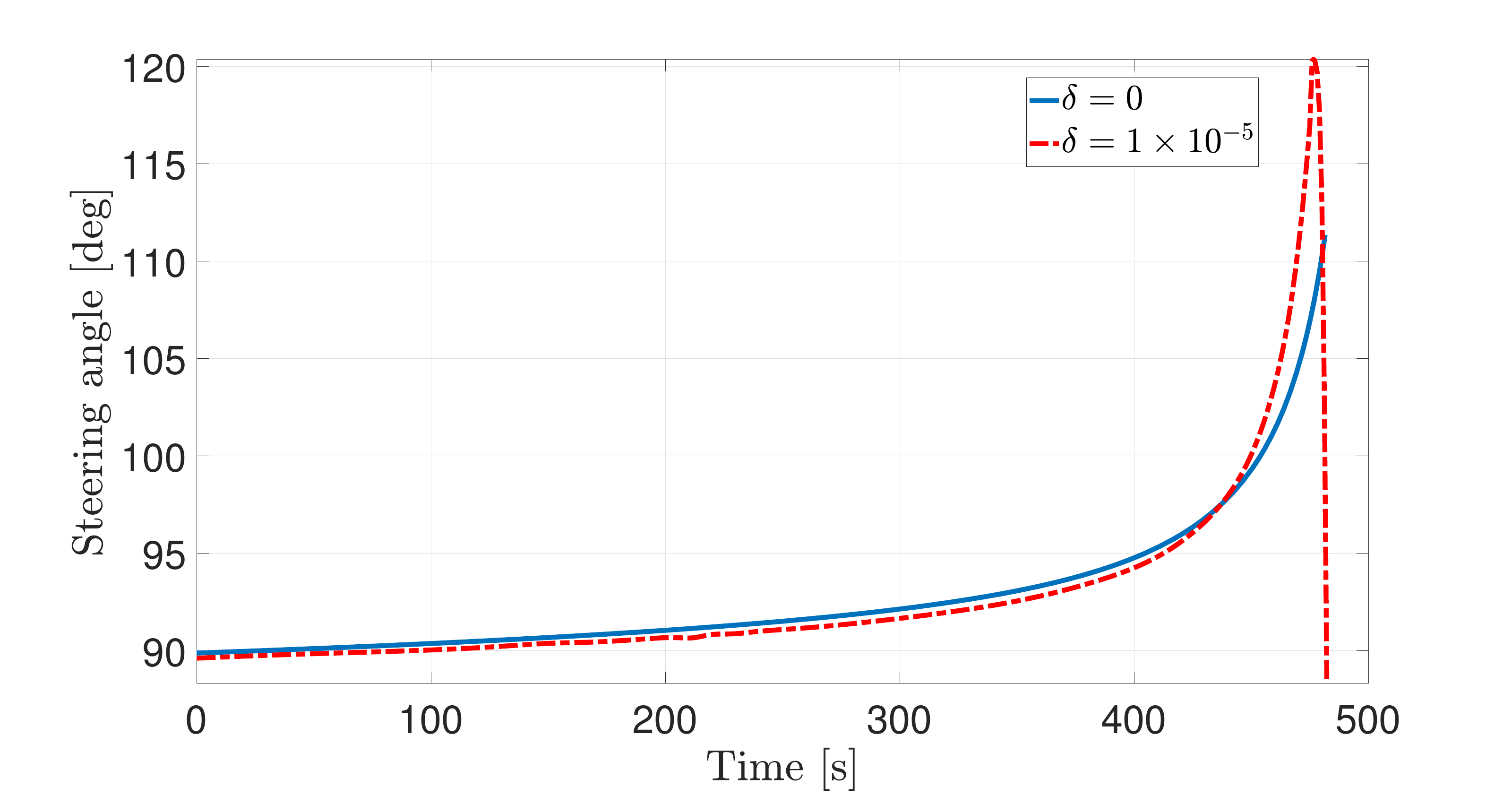}
\caption{Profiles of the steering angle without and with the final steering angle constraint for Case 1.}\label{Fig:case1_control}
\end{center}
\end{figure}

Regarding Case 1, Fig.~\ref{Fig:cooperative_profile} compares the profiles of the altitude, transverse speed, radial speed, and mass. It can be observed that the results without and with the final steering angle constraint are nearly identical, except for the transverse speed.
Fig.~\ref{Fig:case1_control} shows the profiles of the steering angle. It can be seen that the steering angles are quite similar until the final phase of the landing. For $\delta=1\times 10^{-5}$, the steering angle undergoes a rapid change and ends at $\beta^{\mathcal{N}}(t_f) = 88.31$ deg, which is very close to $90$ deg. In contrast, when the final steering angle constraint is ignored, the final steering angle is $111.35$ deg. Moreover, the final state of the lunar lander guided by the NN is: $h^{\mathcal{N}}(t_f) = 5 $ m, $u^{\mathcal{N}}(t_f) = -0.0496 $ m/s, $v^{\mathcal{N}}(t_f) =  -3.8264$ m/s, $m^{\mathcal{N}}(t_f) = 354.5890$ kg. The final mass related to $\delta=0$ is $354.7502$ kg. It should be noted that the final speed can be further reduced for the altitude-to-go of $5$ m, thanks to the retro-propulsion provided by the engine.

\begin{figure}[!htp]
\centering
\begin{subfigure}[t]{4.05cm}
\centering
\includegraphics[width = 4.05cm]{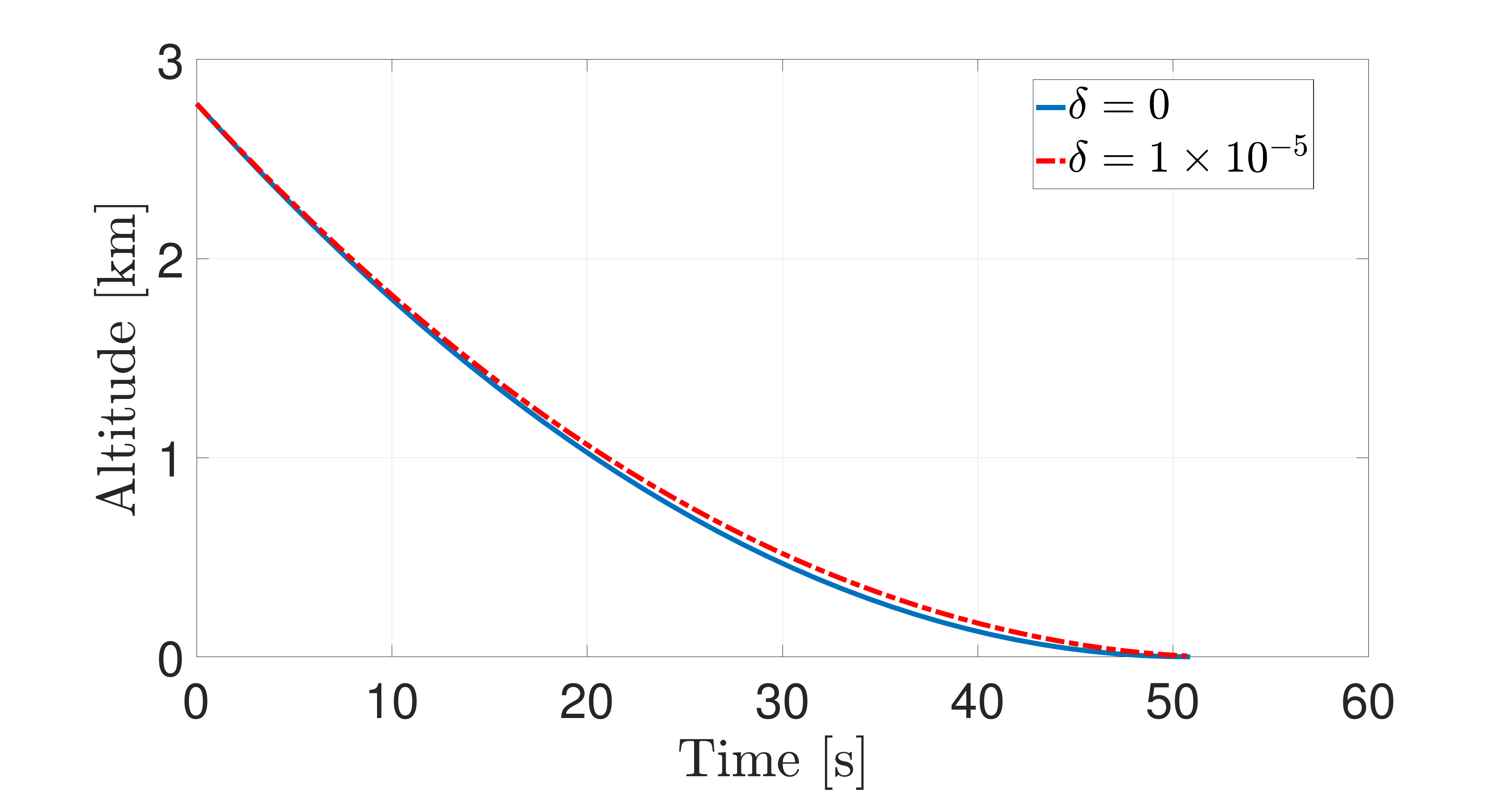}
\caption{Altitude}
\label{Fig:cooperative_control_1_case2}
\end{subfigure}
~~~~~
\begin{subfigure}[t]{4.05cm}
\centering
\includegraphics[width = 4.05cm]{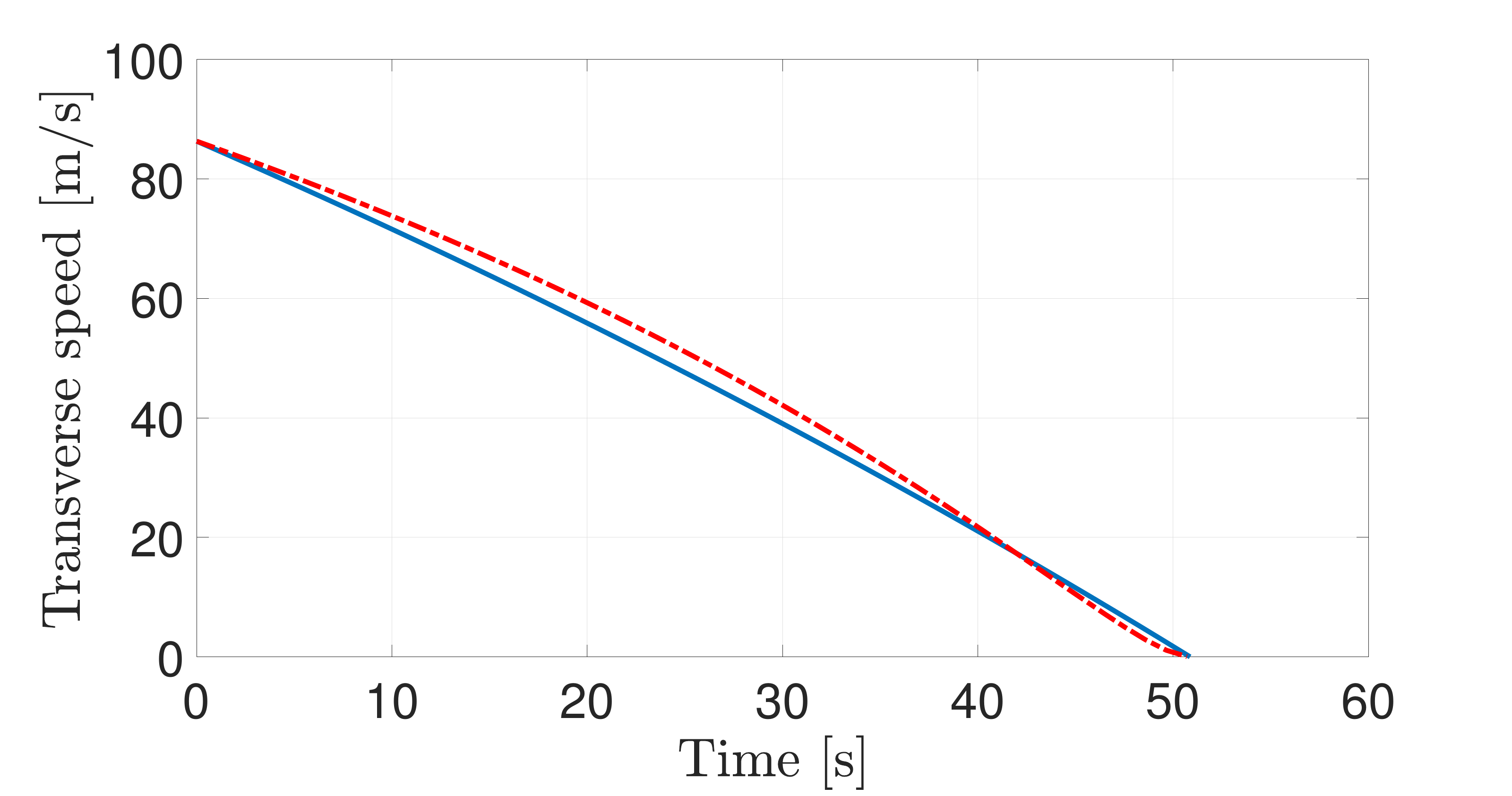}
\caption{Transverse speed}
\label{Fig:cooperative_control_2_case2}
\end{subfigure}\\
\begin{subfigure}[t]{4.05cm}
\centering
\includegraphics[width = 4.05cm]{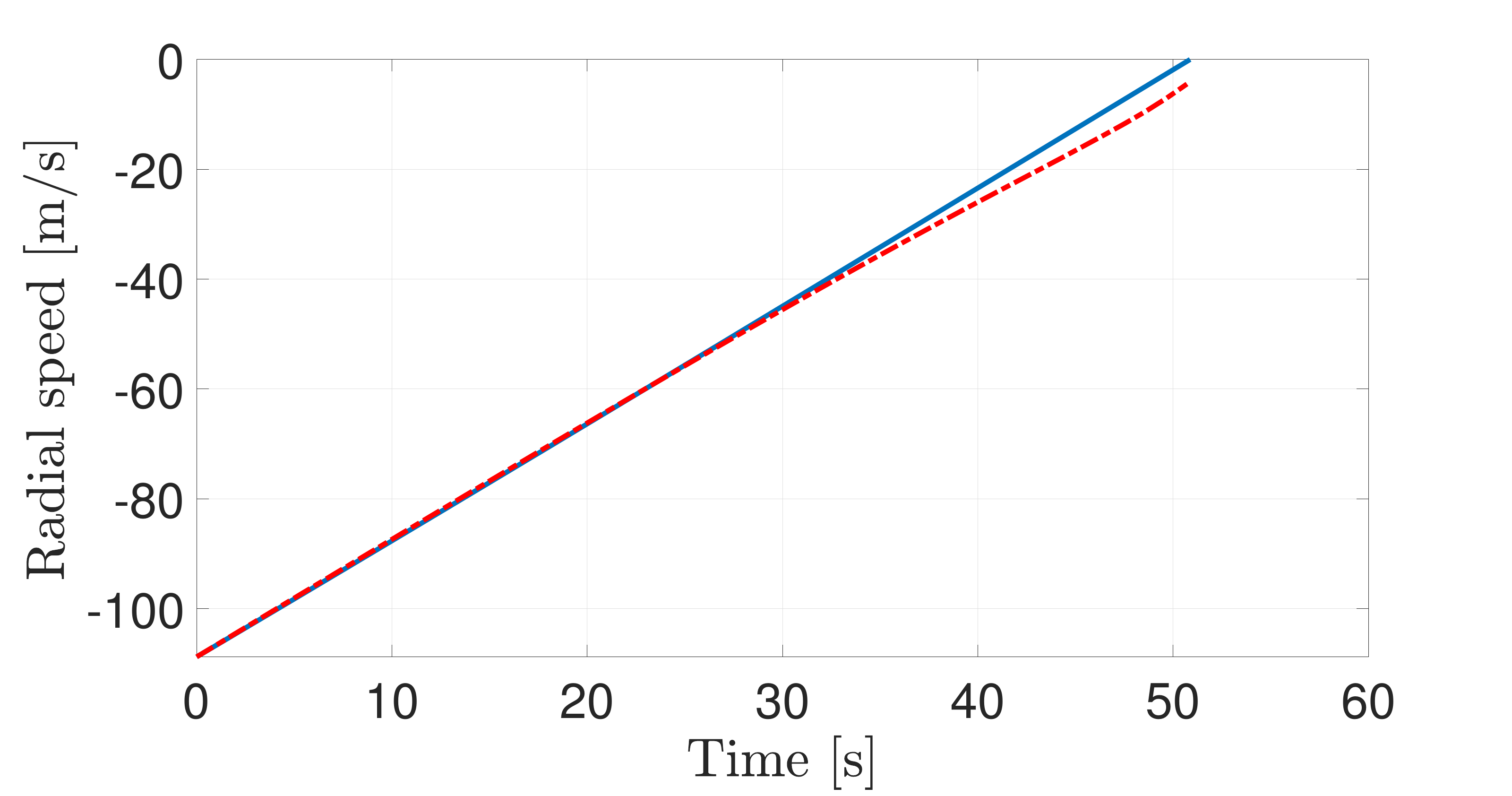}
\caption{Radial speed}
\label{Fig:cooperative_control_3_case2}
\end{subfigure}
~~~~~
\begin{subfigure}[t]{4.05cm}
\centering
\includegraphics[width = 4.05cm]{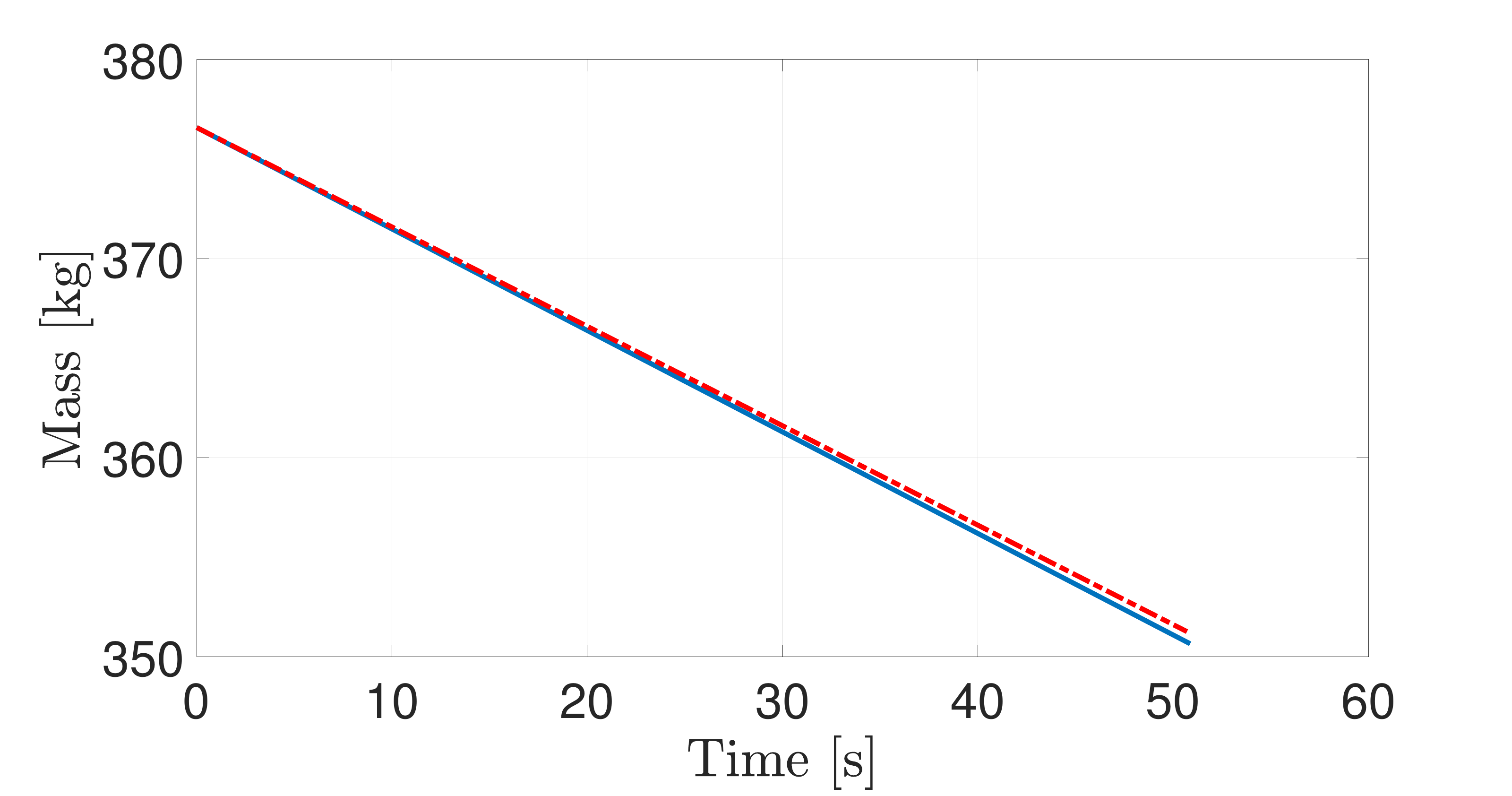}
\caption{Mass}
\label{Fig:cooperative_control_4_case2}
\end{subfigure}
\caption{Profiles of the states without and with the final steering angle constraint for Case 2.}
\label{Fig:cooperative_profile_2}
\end{figure}

\begin{figure}[!htp]
\begin{center}
\includegraphics[scale=0.16]{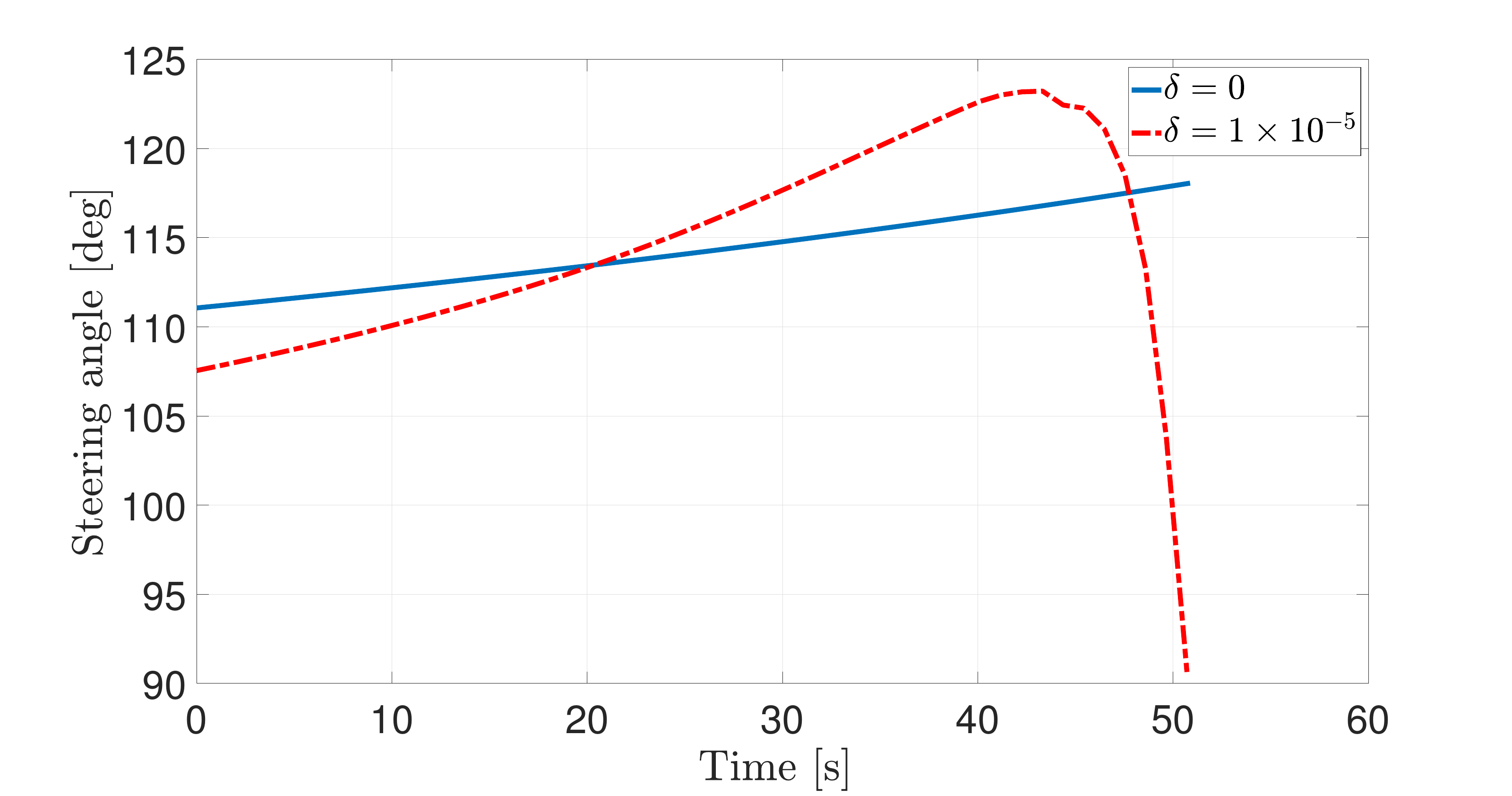}
\caption{Profiles of the steering angle without and with the final steering angle constraint for Case 2.}\label{Fig:case_control_2}
\end{center}
\end{figure}

Fig.~\ref{Fig:cooperative_profile_2} compares the profiles of the altitude, transverse speed, radial speed, and mass for Case 2. 
The corresponding profiles of the steering angle are shown in Fig.~\ref{Fig:case_control_2}. It is evident that the steering angle exhibits a nearly linear progression throughout the entire landing when ignoring the final steering angle constraint. However, when the final steering angle constraint is imposed, the steering angle differs significantly. As a result, the final steering angle for $\delta=1\times 10^{-5}$ is $90.64$ deg. The final state of the lunar lander guided by the NN is: $h^{\mathcal{N}}(t_f) = 5 $ m, $u^{\mathcal{N}}(t_f) = 0.0498 $ m/s, $v^{\mathcal{N}}(t_f) =  -4.4665$ m/s, $m^{\mathcal{N}}(t_f) = 351.267$ kg. Once again, it should be noted that the final speed can be further reduced for the altitude-to-go of $5$ m.
\section{Conclusions}
In this work, in order to achieve a vertical landing, the final attitude constraint was considered as a final control constraint. To address this constraint, an easy-to-implement procedure was proposed to augment the cost functional with a nonnegative small regularization term. Subsequently, a parameterized system was established, allowing for the generation of numerous optimal trajectories. Furthermore, a neural network was trained to approximate the steering angle. Numerical simulations demonstrated that the proposed method was able to generate real-time steering angle commands, guiding the lunar lander to achieve a nearly vertical attitude upon touchdown.
\bibliographystyle{IEEEtran}
\bibliography{IEEEabrv,IEEEexample}

\begin{thebibliography}{10}
\providecommand{\url}[1]{#1}
\csname url@samestyle\endcsname
\providecommand{\newblock}{\relax}
\providecommand{\bibinfo}[2]{#2}
\providecommand{\BIBentrySTDinterwordspacing}{\spaceskip=0pt\relax}
\providecommand{\BIBentryALTinterwordstretchfactor}{4}
\providecommand{\BIBentryALTinterwordspacing}{\spaceskip=\fontdimen2\font plus
\BIBentryALTinterwordstretchfactor\fontdimen3\font minus \fontdimen4\font\relax}
\providecommand{\BIBforeignlanguage}[2]{{%
\expandafter\ifx\csname l@#1\endcsname\relax
\typeout{** WARNING: IEEEtran.bst: No hyphenation pattern has been}%
\typeout{** loaded for the language `#1'. Using the pattern for}%
\typeout{** the default language instead.}%
\else
\language=\csname l@#1\endcsname
\fi
#2}}
\providecommand{\BIBdecl}{\relax}
\BIBdecl

\bibitem{leeghim2016optimal}
H.~Leeghim, D.-H. Cho, and D.~Kim, ``An optimal trajectory design for the lunar vertical landing,'' \emph{Proceedings of the Institution of Mechanical Engineers, Part G: Journal of Aerospace Engineering}, vol. 230, no.~11, pp. 2077--2085, 2016.

\bibitem{lu2018propellant}
P.~Lu, ``Propellant-optimal powered descent guidance,'' \emph{Journal of Guidance, Control, and Dynamics}, vol.~41, no.~4, pp. 813--826, 2018.

\bibitem{cherry1964general}
G.~Cherry, ``A general, explicit, optimizing guidance law for rocket-propelled spaceflight,'' in \emph{Astrodynamics Guidance and Control Conference}, 1964, p. 638.

\bibitem{betts1998survey}
J.~T. Betts, ``Survey of numerical methods for trajectory optimization,'' \emph{Journal of Guidance, Control, and Dynamics}, vol.~21, no.~2, pp. 193--207, 1998.

\bibitem{accikmecse2013lossless}
B.~Acikmese, J.~M. Carson, and L.~Blackmore, ``Lossless convexification of nonconvex control bound and pointing constraints of the soft landing optimal control problem,'' \emph{IEEE Transactions on Control Systems Technology}, vol.~21, no.~6, pp. 2104--2113, 2013.

\bibitem{sagliano2024six}
M.~Sagliano, D.~Seelbinder, S.~Theil, and P.~Lu, ``Six-degree-of-freedom rocket landing optimization via augmented convex--concave decomposition,'' \emph{Journal of Guidance, Control, and Dynamics}, vol.~47, no.~1, pp. 20--35, 2024.

\bibitem{ito2020throttled}
T.~Ito and S.-i. Sakai, ``Throttled explicit guidance to realize pinpoint landing under a bounded thrust magnitude,'' \emph{Journal of Guidance, Control, and Dynamics}, vol.~44, no.~4, pp. 854--861, 2020.

\bibitem{lu2023propellant}
P.~Lu and R.~Callan, ``Propellant-optimal powered descent guidance revisited,'' \emph{Journal of Guidance, Control, and Dynamics}, vol.~46, no.~2, pp. 215--230, 2023.

\bibitem{wang2024new}
K.~Wang, Z.~Chen, Z.~Wei, F.~Lu, and J.~Li, ``A new smoothing technique for bang-bang optimal control problems,'' in \emph{AIAA SCITECH 2024 Forum}, 2024, p. 0727.

\bibitem{sanchez2018real}
C.~S{\'a}nchez-S{\'a}nchez and D.~Izzo, ``Real-time optimal control via deep neural networks: study on landing problems,'' \emph{Journal of Guidance, Control, and Dynamics}, vol.~41, no.~5, pp. 1122--1135, 2018.

\bibitem{wang2022nonlinear}
K.~Wang, Z.~Chen, H.~Wang, J.~Li, and X.~Shao, ``Nonlinear optimal guidance for intercepting stationary targets with impact-time constraints,'' \emph{Journal of Guidance, Control, and Dynamics}, vol.~45, no.~9, pp. 1614--1626, 2022.

\bibitem{WANG2024446}
K.~Wang, F.~Lu, Z.~Chen, and J.~Li, ``Real-time optimal control for attitude-constrained solar sailcrafts via neural networks,'' \emph{Acta Astronautica}, vol. 216, pp. 446--458, 2024.

\bibitem{roozegar2018optimal}
M.~Roozegar, J.~Angeles, and H.~Michalska, ``Optimal control problems with terminal control constraints and benefits of over-actuation,'' in \emph{2018 Annual American Control Conference (ACC)}.\hskip 1em plus 0.5em minus 0.4em\relax IEEE, 2018, pp. 4129--4134.

\bibitem{lee2013polynomial}
C.-H. Lee, T.-H. Kim, M.-J. Tahk, and I.-H. Whang, ``Polynomial guidance laws considering terminal impact angle and acceleration constraints,'' \emph{IEEE Transactions on Aerospace and Electronic Systems}, vol.~49, no.~1, pp. 74--92, 2013.

\bibitem{mcinnes1995path}
C.~R. McInnes, ``Path shaping guidance for terminal lunar descent,'' \emph{Acta Astronautica}, vol.~36, no.~7, pp. 367--377, 1995.

\bibitem{zhou2010optimal}
J.~Zhou, K.~L. Teo, D.~Zhou, and G.~Zhao, ``Optimal guidance for lunar module soft landing,'' \emph{Nonlinear Dynamics and Systems Theory}, vol.~10, no.~2, pp. 189--201, 2010.

\bibitem{sachan2015fuel}
K.~Sachan and R.~Padhi, ``Fuel-optimal {G-MPSP} guidance for powered descent phase of soft lunar landing,'' in \emph{2015 IEEE Conference on Control Applications (CCA)}.\hskip 1em plus 0.5em minus 0.4em\relax IEEE, 2015, pp. 924--929.

\bibitem{pengkun2019}
K.~Peng, R.~Peng, Z.~Huang, and B.~Zhang, ``Implicit shooting method to solve optimal lunar soft landing trajectory,'' \emph{Acta Aeronautica et Astronautica Sinica}, vol.~40, no.~7, pp. 159--167, 2019.

\bibitem{Pontryagin}
L.~S. Pontryagin, V.~G. Boltyanski, R.~V. Gamkrelidze, and E.~F. Mishchenko, \emph{The Mathematical Theory of Optimal Processes (Russian)}.\hskip 1em plus 0.5em minus 0.4em\relax English translation: Interscience, 1962.

\bibitem{zheng2021time}
Y.~Zheng, Z.~Chen, X.~Shao, and W.~Zhao, ``Time-optimal guidance for intercepting moving targets by dubins vehicles,'' \emph{Automatica}, vol. 128, p. 109557, 2021.

\bibitem{chen2019nonlinear}
Z.~Chen and T.~Shima, ``Nonlinear optimal guidance for intercepting a stationary target,'' \emph{Journal of Guidance, Control, and Dynamics}, vol.~42, no.~11, pp. 2418--2431, 2019.

\bibitem{HORNIK1989359}
K.~Hornik, M.~Stinchcombe, and H.~White, ``Multilayer feedforward networks are universal approximators,'' \emph{Neural Networks}, vol.~2, no.~5, pp. 359--366, 1989.

\end{thebibliography}

\end{document}